\documentclass[a4paper,10pt]{article}
\title{Some remarks on representations of Yang-Mills algebras}
\author{Estanislao Herscovich 
\footnote{Departamento de Matem\'atica, FCEyN, Universidad de Buenos Aires, Argentina. 
The author is also a research member of CONICET (Argentina). 
On leave of absence from Institut Joseph Fourier, Universit\'e Grenoble I, France. 
The author would also like to thank the Alexander von Humboldt Foundation and the Bielefeld University, for their support during part of this work. This work was also partially supported by UBACYT 20020130200169BA and PICT 2012-1186.}
}
\date{}

\usepackage{amsmath,amsthm,amsfonts,amssymb,stmaryrd}
\usepackage{latexsym}
\usepackage{color}
\usepackage{graphicx}
\usepackage{float}  
\usepackage{fullpage}
\usepackage[small, bf, margin=90pt, tableposition=bottom]{caption}
\usepackage[T1]{fontenc}
\usepackage{palatino}
\usepackage[fleqn,tbtags]{mathtools}

\usepackage[
  breaklinks,
  naturalnames,
  ]{hyperref}

\usepackage{amsrefs}

\input{xy}
\xyoption{all}
\xyoption{poly}
\usepackage[all]{xy}

\newtheorem{theorem}{Theorem}[section]
\newtheorem{theorem*}{Theorem*}
\newtheorem{corollary}[theorem]{Corollary}

\newtheorem{proposition}[theorem]{Proposition}
\newtheorem{remark}[theorem]{Remark}

\newtheorem{example}[theorem]{Example}

\newtheorem{fact}[theorem]{Fact}
\numberwithin{equation}{section}                    

\let\oldqed\qed
\renewcommand\qed{\oldqed\par\bigskip}

\newcommand\cl[1]{{\langle#1\rangle}}

\newcommand\CC{{\mathbb{C}}}

\newcommand\ZZ{{\mathbb{Z}}}
\newcommand\NN{{\mathbb{N}}}

\newcommand\sll{{\mathfrak{sl}}}

\def\U{{\mathcal U}}

\def\YM{{\mathrm {YM}}}

\def\h{{\mathfrak h}}

\def\ym{{\mathfrak{ym}}}
\def\f{{\mathfrak f}}

\def\MyNode{\ifcase\xypolynode\or
      (W \otimes X) \otimes (Y \otimes Z)
    \or
      ((W \otimes X) \otimes Y) \otimes Z
    \or
      (W \otimes (X \otimes Y)) \otimes Z
    \or
      W \otimes ((X \otimes Y) \otimes Z)
    \or
      W \otimes (X \otimes (Y \otimes Z))
    \fi
  }%

\def\MyNodes{\ifcase\xypolynode
    \or
      X \otimes Y
    \or
      X \otimes (e \otimes Y)
    \or
     (X \otimes e) \otimes Y
    \fi
  }%

\def\MyNodessr{\ifcase\xypolynode
    \or
      X \otimes Y
    \or
      e \otimes (X \otimes Y)
    \or
     (e \otimes X) \otimes Y
    \fi
  }%

\def\MyNodessl{\ifcase\xypolynode
    \or
      X \otimes Y
    \or
      (X \otimes Y) \otimes e
    \or
       X \otimes (Y \otimes e)
    \fi
  }%
\begin{document}

\maketitle
\begin{abstract}
   {
In this article we present some probably unexpected (in our opinion) properties of representations of Yang-Mills algebras. 
We first show that any free Lie algebra with $m$ generators is a quotient of the Yang-Mills algebra $\ym(n)$ on $n$ generators, for $n \geq 2 m$. 
We derive from this that any semisimple Lie algebra, and even any affine Kac-Moody algebra is a quotient 
of $\ym(n)$, for $n \geq 4$. 
Combining this with previous results on representations of Yang-Mills algebras given in \cite{HS10}, one may obtain solutions to the Yang-Mills equations by differential operators acting on sections of twisted vector bundles on the affine space of dimension $n \geq 4$ associated to representations of any semisimple Lie algebra. 
We also show that this quotient property does not hold for $n = 3$, since any morphism of Lie algebras from $\ym(3)$ to 
$\sll(2,k)$ has in fact solvable image. 
   }
\end{abstract}



\section{Introduction}

The purpose of this short article is to show some unobserved properties about the representation theory of Yang-Mills algebras.
In the previous article \cite{HS10}, we have already considered the question of finding representations of Yang-Mills algebras, and we have obtained several of them factoring through Weyl algebras. 
Here we proceed in a more simple fashion. 
It aims to answer a question posed by J. Alev of finding which (semisimple) Lie algebras can be obtained as quotients of Yang-Mills algebras $\ym(n)$, for $n \geq 3$. 
The case $n = 2$ is already well-known, for the Yang-Mills algebra $\ym(2)$ is isomorphic to the first Heisenberg Lie algebra. 
Even though some of the results are somehow elementary, the consequences we derive from them seem unexpected (see for instance Remark \ref{remark:repimp}). 
In particular one observes a stark difference between the cases $n = 3$ and $n \geq 4$ (see Proposition \ref{proposition:ym3}), together with a rather particular behaviour in the former case. 

We recall that Yang-Mills algebras we introduced by A. Connes and M. Dubois-Violette in \cite{CD02}. 
As explained there (see also the introduction of \cite{HS10}) the interest in them rely on the fact that one may provide a ``(local) background independent'' description of gauge theory in physics. 
Indeed, the relations defining the Yang-Mills algebras on $n$ generators are just generated by the Euler-Lagrange equations satisfied by the components of the covariant derivative on a complex vector bundle 
over the $n$-dimensional affine space. 
Moreover, as shown by N. Nekrasov they also arise as limits of algebras appearing in the gauge theory of $D$-branes and open string theory in \cite{Ne03}. 
In that article the mentioned author stated the little knowledge about the representation theory of these algebras. 

I would like to express my deep gratitude to Jacques Alev, for introducing me to this field, and sharing his ideas and profound understanding of Lie algebras. 
Moreover, I would also like to deeply thank Rupert Yu, for several discussions and comments, and for his interesting example that is presented at the end of this paper.  

\section{Yang-Mills algebras}

Throughout this article $k$ will denote an algebraically closed field of characteristic zero.
Let $n$ be a positive integer such that $n \geq 2$ and let $\f(x_{1}, \dots, x_{n})$ be the free Lie algebra on $n$ generators $\{ x_{1}, \dots, x_{n} \}$.
This Lie algebra is provided with a locally finite dimensional $\NN$-grading.

Following \cite{CD02}, the quotient Lie algebra
\[     \ym(n) = \f(x_{1}, \dots, x_{n})/\cl{\{ \sum_{i=1}^{n} [x_{i},[x_{i},x_{j}]] : 1 \leq j \leq n \}}     \]
is called the \emph{Yang-Mills algebra on $n$ generators}.

The $\NN$-grading of $\f(x_{1}, \dots, x_{n})$ induces an $\NN$-grading of $\ym(n)$, which is also locally finite dimensional.
We denote $\ym(n)_{j}$ the $j$-th homogeneous component and
\begin{equation}
\label{eq:gradym}
\ym(n)^{l} = \bigoplus_{j = 1}^{l} \ym(n)_{j}.
\end{equation}

The enveloping algebra $\U(\ym(n))$ will be denoted $\YM(n)$.
It is the \emph{(associative) Yang-Mills algebra on $n$ generators}.
If $V(n) = \text{span}_{k} \cl{\{ x_{1} , \dots , x_{n} \}}$, we see that
\[     \YM(n) \simeq T(V(n))/\cl{R(n)},     \]
where $T(V(n))$ denotes the tensor $k$-algebra on $V(n))$ and 
\begin{equation}
\label{eq:relym}
     R(n) = \mathrm{span}_{k} \cl{\{ \sum_{i=1}^{n} [x_{i},[x_{i},x_{j}]] : 1 \leq j \leq n \}} \subseteq V(n)^{\otimes 3}.
\end{equation}
We shall denote by $r_{j} = \sum_{i=1}^{n} [x_{i},[x_{i},x_{j}]]$, for $j = 1, \dots, n$, the $r$-th basis element of space of relations $R(n)$. 

We also recall that the obvious projection $V(n) \rightarrow V(m)$, for $m \leq n$, given by $x_{i} \mapsto x_{i}$, for $i = 1, \dots, m$, and 
$x_{i} \mapsto 0$, for $i = m+1, \dots, n$, induces a surjective morphism of Lie algebras $\ym(n) \rightarrow \ym(m)$, 
and also a surjective morphism of algebras $\YM(n) \rightarrow \YM(m)$. 

As noted in \cite{HS10}, the algebra $\ym(n)$ is nilpotent if $n=2$, in which case it is also finite dimensional
(see \cite{HS10}, Example 2.1). 
In fact, it is isomorphic to the first Heisenberg Lie algebra $\h_{1}$.
When $n \geq 3$, $\ym(n)$ is neither finite dimensional nor nilpotent (see \cite{HS10}, Remark 3.14).
We shall see however that there is an important difference between the case $n = 3$ and $n \geq 4$. 

\begin{proposition}
Let $m \in \NN$ and $n \geq 2 m$. 
Then there exists a surjective morphism of (graded) Lie algebras from $\ym(n)$ to $\f(m)$, the free Lie algebra on $m$ generators. 
\end{proposition}
\noindent\textbf{Proof.}
Since $\ym(2m)$ is an epimorphic image of the graded Lie algebra $\ym(n)$, for $n \geq 2m$, it suffices to find a surjective morphism of graded Lie algebras
$\ym(2m) \rightarrow \f(m)$. 
Let us denote the generators of $\f(m)$ by $y_{1}, \dots, y_{m}$. 
The linear map from $V(n)$ to $\f(m)$ given by 
\begin{align*}
 x_{j} &\mapsto y_{j},
 \\
 x_{m+j} &\mapsto i y_{j},
\end{align*}
for $j = 1, \dots, m$, induces a surjective morphism of graded Lie algebras $\phi : \f(x_{1}, \dots, x_{2 m}) \rightarrow \f(m)$. 
It obviously satisfies that $\phi(R(2 m)) = 0$, since 
\[     \phi\Big(\sum_{j=1}^{2 m} [x_{j},[x_{j},x_{l}]]\Big) = \sum_{j=1}^{2 m} [\phi(x_{j}),[\phi(x_{j}),\phi(x_{l})]] 
           = \sum_{j=1}^{m} [y_{j},[y_{j},\phi(x_{l})]] - \sum_{j=m+1}^{2 m} [y_{j-m},[y_{j-m},\phi(x_{l})]] = 0,     \]
for all $l = 1, \dots, 2 m$. 
Hence, we get a surjective morphism of graded Lie algebras $\phi : \ym(2 m) \rightarrow \f(m)$, and the proposition follows. 
\qed

We recall that the \emph{Witt algebra} $\mathfrak{W}$ is the Lie algebra over $k$ of derivations of $k[z,z^{-1}]$, 
(which, if $k = \CC$, is just the complexification of the Lie algebra of polynomial vector fields on the circle). 
It was introduced by E. Cartan in \cite{Ca1909}, 
and it can be explicitly described as the vector space over $k$ with basis $\{ e_{n} \}_{n \in \ZZ}$ and 
brackets $[e_{n} , e_{m}] = (m-n) e_{m+n}$, for all $m, n \in \ZZ$. 
The central extension of $\mathfrak{W}$ given by $\mathrm{Vir} = \mathfrak{W} \oplus k.c$ with the bracket 
\[     [e_{n},e_{m}] = (m-n) e_{m+n} + \frac{\delta_{m+n,0} (m^{3}-m)}{12}  c,  \hskip 1cm [e_{n},c] = 0,     \]
for $m, n \in \ZZ$ is called the \emph{Virasoro algebra}, and it was introduced by I. Gelfand and D. Fuks in \cite{GF68}. 
It is easy to see that both Lie algebras may be generated by two elements, \textit{e.g.} $e_{-2}$ and $e_{3}$. 
As a consequence of the previous proposition we get the following result. 
\begin{corollary}
If $n \geq 4$, both the Witt algebra $\mathfrak{W}$ and the Virasoro alegbra $\mathrm{Vir}$ are quotients of the Yang-Mills algebra $\ym(n)$.
\end{corollary}

As another consequence of the previous proposition and the fact proved by M. Kuranishi that any semisimple Lie algebra over a field of characteristic zero 
is generated by two elements (see \cite{Ku51}, Thm. 6), we also obtain the following result.
\begin{corollary}
\label{corollary:ss}
Let $n \geq 4$. 
Then, any semisimple Lie algebra is a quotient of the Yang-Mills algebra $\ym(n)$.
\end{corollary}

Furthermore, let us consider $A = (a_{ij})_{i,j = 1, \dots, m}$ be a matrix of entries in $k$ of rank $r$, and let 
$(\mathfrak{h}, \Pi, \Pi^{\vee})$ be a \emph{realization of $A$}, \textit{i.e.} $\mathfrak{h}$ is a vector space, $\Pi =\{ \alpha_{1}, \dots, \alpha_{m} \} \subseteq \mathfrak{h}^{*}$ and $\Pi^{\vee} =\{ \alpha_{1}^{\vee}, \dots, \alpha_{m}^{\vee} \} \subseteq \mathfrak{h}$ are linearly independent subsets of the dual $\mathfrak{h}^{*}$ of $\mathfrak{h}$ and of $\mathfrak{h}$, respectively, satisfying that 
\[     \langle \alpha_{i}^{\vee}, \alpha_{j} \rangle = a_{ij}, \text{ for all $i,j = 1, \dots, m$}, \hskip 1cm \text{and} \hskip 1cm
\dim \mathfrak{h} = 2 m - r.     \] 
This notion was introduced by V. Kac in \cite{Ka}, Ch. 1. 
Moreover, the definition of \emph{morphism of realizations of a fixed matrix} is straightforward. 
In fact, one can prove that the realization of a matrix is unique up (not necessarily unique) isomorphism (see \cite{Ka}, Prop. 1.1). 
One defines the Lie algebra $\tilde{\mathfrak{g}}(A)$ over $k$ generated by the vector space 
given by the direct sum of  $\mathfrak{h}$ and the one spanned by $\{ e_{i}, f_{i} \}_{i = 1, \dots, m}$, 
modulo the relations 
\begin{equation*}
\begin{split}
     [e_{i},f_{j}] &= \delta_{i,j} \alpha_{i}^{\vee},    \hskip 1cm [h,h'] = 0, 
     \\
     [h,e_{i}] &= \langle h , \alpha_{i} \rangle e_{i}, \hskip 0.7cm [h,f_{i}] = - \langle h , \alpha_{i} \rangle f_{i},
\end{split}
\end{equation*} 
for all $i,j = 1, \dots, m$ and $h, h' \in \mathfrak{h}$. 
One may prove that the canonical map from $\mathfrak{h}$ to $\tilde{\mathfrak{g}}(A)$ is injective (see \cite{Ka}, Thm. 1.2, (a)), so we may regard $\mathfrak{h}$ inside of $\tilde{\mathfrak{g}}(A)$. 
There exists a unique maximal element $\mathfrak{r}$ in the set of ideals of the Lie algebra $\tilde{\mathfrak{g}}(A)$ intersecting $\mathfrak{h}$ trivially (see \cite{Ka}, Thm. 1.2, (e)), and define $\mathfrak{g}(A) = \tilde{\mathfrak{g}}(A)/\mathfrak{r}$. 
It is called the \emph{Kac-Moody algebra associated to the matrix $A$} in case $A$ is a \emph{generalized Cartan matrix}, \textit{i.e.} if $A$ has integer entries, $a_{ii} = 2$, for all $i = 1, \dots, m$, $a_{ij} \leq 0$, for all $i \neq j$, and $a_{ij} = 0$ implies $a_{ji} = 0$. 
By extending a theorem of J.-P. Serre, it is known that finite dimensional semisimple Lie algebras are exactly Kac-Moody algebras for a nonsingular (generalized Cartan) matrix (see \cite{Ka}, Prop. 4.9). 
Thus, the previous corollary can be strengthen by using the theorem in \cite{LW86} (see also the last remark of that article). 
\begin{corollary}
\label{corollary:km}
Let $n \geq 4$. 
Consider a matrix $A = (a_{ij})_{i,j = 1, \dots, m}$ of entries in $k$ of rank $r$. 
Then, any of the Lie algebras $\mathfrak{g}(A)$ or $\tilde{\mathfrak{g}}(A)$ are a quotient of the Yang-Mills algebra $\ym(n)$ if $r + 2 \geq m$ (this includes in particular all affine Kac-Moody algebras). 
On the other hand, if $r + 2 \leq m$, then either $\mathfrak{g}(A)$ or $\tilde{\mathfrak{g}}(A)$ is a quotient of the Yang-Mills algebra $\ym(n)$ for $n \geq 2(m-r)$. 
\end{corollary}

\begin{remark}
\label{remark:repimp}
By combining the result of the previous corollary together with \cite{HS10}, Thm 1.1, we obtain 
solutions of the Yang-Mills equations by differential operators acting on sections of twisted vector bundles on the affine space of dimension $n \geq 4$ associated to representations of any semisimple Lie algebra (or, also, of any affine Kac-Moody algebra).    
\end{remark}

We would like to remark that the last two corollaries are not true for $n = 3$, as we now show. 
In order to prove that, we shall make use of the following simple result.
\begin{fact}
\label{fact:1}
Let us consider the vector space $k^{2}$ provided with the canonical nondegenerate symmetric bilinear form $x \cdot y = x_{1} y_{1} + x_{2} y_{2}$, for 
$x=(x_{1},x_{2})$ and $y=(y_{1},y_{2})$. 
Let $x, y \in k^{2}$ be two elements with $x$ nonzero, such that $x \cdot x = 0$ and $x \cdot y = 0$. 
Then $y$ is a scalar multiple of $x$, and in particular $y \cdot y = 0$. 
\end{fact}

We have now the:
\begin{proposition}
\label{proposition:ym3}
Let $\phi : \ym(3) \rightarrow \sll(2,k)$ be any morphism of Lie algebras. 
Then, the image of $\phi$ is a solvable Lie subalgebra of $\sll(2,k)$. 
\end{proposition}
\noindent\textbf{Proof.}
We recall that $\sll(2,k)$ is the vector space spanned by $\{e, h, f\}$ with relations $[h,e] = 2 e$, $[h,f] = -2 f$, and $[e,f]=h$. 
If $\dim_{k}(\phi(V(3))) = 0$, there is nothing to prove, because $\mathrm{Im}(\phi)$ is the subalgebra of $\sll(2,k)$ generated by $\phi(V(3))$, 
which in this case is even nilpotent. 

Let us thus suppose that $\dim_{k}(\phi(V(3))) \neq 0$. 
Then there exists $i \in \{ 1, 2, 3 \}$ such that $\phi(x_{i}) \neq 0$. 
Without loss of generality, we may assume that $i = 3$. 
Since any nonzero element $z$ of $\sll(2,k)$ is either nilpotent or semisimple (\textit{i.e.} the map $\mathrm{ad}(z) : \sll(2,k) \rightarrow \sll(2,k)$ is nilpotent or diagonalizable, respectively), we will consider two separated cases. 

First, let us suppose that $\phi(x_{3})$ is nilpotent.   
By the Jacobson-Morozov theorem (see \cite{Kn02}, Thm. 10.3), we may even suppose that $\phi(x_{3}) = e$. 
Set $\phi(x_{i}) = \alpha_{i} e + \beta_{i} h + \gamma_{i} f$, for $i = 1, 2$. 
Let us denote $\alpha = (\alpha_{1},\alpha_{2})$, $\beta = (\beta_{1},\beta_{2})$, $\gamma = (\gamma_{1},\gamma_{2})$ the corresponding vectors in $k^{2}$, 
and we shall denote by $x \cdot y = x_{1} y_{1} + x_{2} y_{2}$ the canonical nondegenerate symmetric bilinear form of $k^{2}$, as in Fact \ref{fact:1}. 
In this case we may further assume that $\gamma \neq 0$, because, if not, we have that $\mathrm{Im}(\phi) \subseteq k.e \oplus k.h$, which is solvable. 
We shall prove that the nonvanishing of $\gamma$ gives an absurd. 

On the one hand, we have that 
\begin{align*}
   \phi\Big(\sum_{i=1}^{3} [x_{i},[x_{i},x_{3}]]\Big) &= \sum_{i=1}^{2} [\alpha_{i} e + \beta_{i} h + \gamma_{i} f,[\alpha_{i} e + \beta_{i} h + \gamma_{i} f,e]]
   = \sum_{i=1}^{2} [\alpha_{i} e + \beta_{i} h + \gamma_{i} f, 2 \beta_{i} e - \gamma_{i} h] 
   \\
   &= \sum_{i=1}^{2} \big(2 (\alpha_{i} \gamma_{i} + 2 \beta_{i}^{2}) e - \beta_{i} \gamma_{i} h - 2 \gamma_{i}^{2} f\big), 
\end{align*}
which implies that $\phi(r_{3}) = 0$ if and only if 
\begin{equation}
\label{eq:rel3nil}
   2 \beta \cdot \beta + \alpha \cdot \gamma = 0, \hskip 0.5cm \beta \cdot \gamma = 0, \hskip 0.5cm \gamma \cdot \gamma = 0.   
\end{equation}
By the Fact \ref{fact:1} together with the previous two last equalities, we conclude that $\beta$ is a multiple of $\gamma$, and thus $\beta \cdot \beta = 0$. 
By applying this last equality to the first of the identities in \eqref{eq:rel3nil}, we further conclude that $\alpha \cdot \gamma = 0$, which in turn implies that $\alpha$ is also a multiple of $\gamma$ and satisfies $\alpha \cdot \alpha = 0$.

On the other hand, if $j \in \{ 1, 2 \}$, we obtain 
\begin{equation}
\label{eq:cuentanil}
\begin{split}
   \phi\Big(&\sum_{i=1}^{3} [x_{i},[x_{i},x_{j}]]\Big) 
   \\
  &= \sum_{i=1}^{2} [\alpha_{i} e + \beta_{i} h + \gamma_{i} f,[\alpha_{i} e + \beta_{i} h + \gamma_{i} f,\alpha_{j} e + \beta_{j} h + \gamma_{j} f]] + [e,[e,\alpha_{j} e + \beta_{j} h + \gamma_{j} f]]
   \\
   &=  \sum_{i=1}^{2} \Big([\alpha_{i} e + \beta_{i} h + \gamma_{i} f, 2 (\beta_{i} \alpha_{j} - \beta_{j} \alpha_{i}) e + (\alpha_{i} \gamma_{j} - \alpha_{j} \gamma_{i}) h + 2 (\gamma_{i} \beta_{j} - \gamma_{j} \beta_{i})f]\Big) - 2 \gamma_{j} e
   \\
   &= \sum_{i=1}^{2} \Big(2 \big(2 \beta_{i} (\beta_{i} \alpha_{j} - \beta_{j} \alpha_{i}) - \alpha_{i} (\alpha_{i} \gamma_{j} - \alpha_{j} \gamma_{i})\big) e + 2 \big(\alpha_{i} (\gamma_{i} \beta_{j} - \gamma_{j} \beta_{i}) - \gamma_{i} (\beta_{i} \alpha_{j} - \beta_{j} \alpha_{i}) \big) h 
  \\
  &\phantom{= \sum_{i=1}^{2}} 
+ 2 \big(\gamma_{i} (\alpha_{i} \gamma_{j} - \alpha_{j} \gamma_{i}) - 2 \beta_{i} (\gamma_{i} \beta_{j} - \gamma_{j} \beta_{i})\big) f \Big) - 2 \gamma_{j} e, 
\end{split}
\end{equation}
so $\phi(r_{j})$ vanishes for $j = 1, 2$ if and only if 
\begin{equation}
\label{eq:reljnil*}
\begin{split}
   \big(2 (\beta \cdot \beta) + (\alpha \cdot \gamma) \big) \alpha - 2 (\beta \cdot \alpha) \beta - \big((\alpha \cdot \alpha) + 1\big) \gamma &= 0, 
   \\
   2 (\alpha \cdot \gamma) \beta - (\alpha \cdot \beta) \gamma - (\beta \cdot \gamma) \alpha &= 0, 
   \\
   \big(2 (\beta \cdot \beta) + (\alpha \cdot \gamma) \big) \gamma - (\gamma \cdot \gamma) \alpha - 2 (\beta \cdot \gamma) \beta &= 0.   
\end{split}
\end{equation}
Using Eq. \eqref{eq:rel3nil}, together with the comments after it, we see that the previous collection 
of identities is equivalent to $\gamma = 0$, which is absurd. 
The contradiction comes from the assumption that $\gamma \neq 0$, so we see that if $\phi(x_{3})$ is nilpotent, then $\mathrm{Im}(\phi)$ is solvable. 

Finally, let us suppose that $\phi(x_{3})$ is semisimple. 
Since all Cartan subalgebras of a simple Lie algebra are conjugated by an inner automorphism (see \cite{Kn02}, Thm. 2.15), we may assume that, after applying such an isomorphism, 
$\phi(x_{3}) = h$. 
As before, we write $\phi(x_{i}) = \alpha_{i} e + \beta_{i} h + \gamma_{i} f$, for $i = 1, 2$, and 
denote $\alpha = (\alpha_{1},\alpha_{2})$, $\beta = (\beta_{1},\beta_{2})$, $\gamma = (\gamma_{1},\gamma_{2})$ the corresponding vectors in $k^{2}$.  

By a similar computation as in the previous situation, we have that 
\begin{align*}
   \phi\Big(\sum_{i=1}^{3} [x_{i},[x_{i},x_{3}]]\Big) &= \sum_{i=1}^{2} [\alpha_{i} e + \beta_{i} h + \gamma_{i} f,[\alpha_{i} e + \beta_{i} h + \gamma_{i} f,h]]
   = \sum_{i=1}^{2} [\alpha_{i} e + \beta_{i} h + \gamma_{i} f, - 2 \alpha_{i} e + 2 \gamma_{i} f] 
   \\
   &= 4 \sum_{i=1}^{2} \Big(-(\alpha_{i} \beta_{i}) e + (\alpha_{i} \gamma_{i}) h - (\beta_{i} \gamma_{i}) f \Big), 
\end{align*}
which implies that $\phi(r_{3}) = 0$ if and only if 
\begin{equation}
\label{eq:rel3semi}
   \alpha \cdot \beta = 0, \hskip 0.5cm \alpha \cdot \gamma = 0, \hskip 0.5cm \beta \cdot \gamma = 0.   
\end{equation}
On the other hand, if $j \in \{ 1, 2 \}$, 
\begin{equation}
\label{eq:cuentasemi}
\begin{split}
   \phi\Big(&\sum_{i=1}^{3} [x_{i},[x_{i},x_{j}]]\Big) 
   \\
   &= \sum_{i=1}^{2} [\alpha_{i} e + \beta_{i} h + \gamma_{i} f,[\alpha_{i} e + \beta_{i} h + \gamma_{i} f,\alpha_{j} e + \beta_{j} h + \gamma_{j} f]] + [h,[h,\alpha_{j} e + \beta_{j} h + \gamma_{j} f]]
   \\
   &= \sum_{i=1}^{2} [\alpha_{i} e + \beta_{i} h + \gamma_{i} f, 2 (\beta_{i} \alpha_{j} - \beta_{j} \alpha_{i}) e + (\alpha_{i} \gamma_{j} - \alpha_{j} \gamma_{i}) h + 2 (\gamma_{i} \beta_{j} - \gamma_{j} \beta_{i})f] 
  + 4 \alpha_{j} e + 4 \gamma_{j} f,
   \\
   &= \sum_{i=1}^{2} \Big(2 \big(2 \beta_{i} (\beta_{i} \alpha_{j} - \beta_{j} \alpha_{i}) - \alpha_{i} (\alpha_{i} \gamma_{j} - \alpha_{j} \gamma_{i})\big) e 
 + 2 \big(\alpha_{i} (\gamma_{i} \beta_{j} - \gamma_{j} \beta_{i}) - \gamma_{i} (\beta_{i} \alpha_{j} - \beta_{j} \alpha_{i}) \big) h 
   \\
   &\phantom{= \sum_{i=1}^{2}} + 2 (\big(\gamma_{i} (\alpha_{i} \gamma_{j} - \alpha_{j} \gamma_{i}) - 2 \beta_{i} (\gamma_{i} \beta_{j} - \gamma_{j} \beta_{i})\big) f\Big) + 4 \alpha_{j} e + 4 \gamma_{j} f. 
\end{split}
\end{equation}
So, $\phi(r_{j}) = 0$, for $j = 1, 2$, becomes in this case equivalent to the following collection of identities
\begin{equation}
\label{eq:reljsemi*}
\begin{split}
   \big(2 (\beta \cdot \beta) + 2 + (\alpha \cdot \gamma) \big) \alpha - 2 (\beta \cdot \alpha) \beta - (\alpha \cdot \alpha) \gamma &= 0, 
   \\
   2 (\alpha \cdot \gamma) \beta - (\alpha \cdot \beta) \gamma - (\beta \cdot \gamma) \alpha &= 0, 
   \\
   \big(2 (\beta \cdot \beta) + 2 + (\alpha \cdot \gamma) \big) \gamma - (\gamma \cdot \gamma) \alpha - 2 (\beta \cdot \gamma) \beta &= 0.   
\end{split}
\end{equation}
By application of \eqref{eq:rel3semi}, they can be reduced to 
\begin{equation}
\label{eq:reljsemi}
\begin{split}
   2 \big(1 + (\beta \cdot \beta)\big) \alpha &= (\alpha \cdot \alpha) \gamma, 
   \\
   2 \big(1 + (\beta \cdot \beta)\big) \gamma &= (\gamma \cdot \gamma) \alpha. 
\end{split} 
\end{equation}
Now, using that $\alpha \cdot \gamma = 0$, any of the previous equations implies that $(\alpha \cdot \alpha) (\gamma \cdot \gamma) = 0$, so 
either $\alpha \cdot \alpha = 0$ or $\gamma \cdot \gamma = 0$. 
In any case, Fact \ref{fact:1} implies that both $\alpha \cdot \alpha = 0$ and $\gamma \cdot \gamma = 0$ hold. 
Moreover, since $\alpha \cdot \beta = 0$, another application of Fact \ref{fact:1} implies that $\beta \cdot \beta = 0$. 
Hence, \eqref{eq:reljsemi} gives $\alpha = \gamma = 0$. 
This is an absurd, for we have supposed that $\gamma \neq 0$. 
The proposition is thus proved. 
\qed

\begin{remark}
The morphism from $\ym(3)$ to $\sll(2,k)$ given by $x_{1} \mapsto h$, $x_{2} \mapsto e$, $x_{3} \mapsto i h$ has solvable image, which is not nilpotent. 
\end{remark}

Interestingly, we have however a surjective morphism from $\ym(3)$ to $\sll(3,k)$, which we now explain. 
It was observed by R. Yu and we are indebted to him for it. 
\begin{example}
Let us consider the Lie algebra $\tilde{\ym}(3)$ defined as the quotient 
\[     \tilde{\ym}(3) = \f(x_{1}, x_{2}, x_{3})/\cl{\{ [x_{i},[x_{i},x_{j}]] : 1 \leq i, j \leq n \}}.     \]
It is clearly a quotient of $\ym(3)$. 
We shall show that $\sll(3,k)$ is a quotient of $\tilde{\ym}(3)$. 

As usual, for $i, j = 1, 2, 3$, define $E^{ij}$ the matrix satisfying that its entry $E^{ij}_{mn} = \delta_{i,m} \delta_{j,n}$, for $m, n = 1, 2, 3$. 
Note in particular that $[E^{1 2}, E^{2 3}] = E^{1 3}$, these elements $E^{1 2}$, $E^{2 3}$ and $E^{3 1}$ belong to the Lie algebra $\sll(3,k)$, and in fact they generate the latter. 

Define a $k$-linear map $\psi$ from $V(3)$ to $\sll(3,k)$ by $x_{1} \mapsto E^{1 2}$, $x_{2} \mapsto E^{2 3}$, 
and $x_{3} \mapsto E^{3 1}$. 
By a direct computation it induces a morphism of Lie algebras from $\tilde{\ym}(3)$ to $\sll(3,k)$.  
Since the image of $\psi$ generates $\sll(3,k)$, we get that the latter is a quotient of the Lie algebra $\tilde{\ym}(3)$, 
which in turn implies that $\sll(3,k)$ is also a quotient of the Yang-Mills algebra $\ym(3)$.
\end{example}




\end{document}